\documentclass[10pt,twoside]{article}
\usepackage{graphicx}
\usepackage{Latex-document}

\usepackage{float}
\restylefloat{figure}

\def\be{\begin{equation}}
\def\ee{\end{equation}}
\def\part{\partial}
\def\R{{\mathbf R}}

\def\nuhalf{\nu+1/2}
\def\numhalf{\nu-1/2}
\def\xnu{x_\nu}
\def\xnuhalf{x_{\nuhalf}}

\def\tn{t^n}
\def\tnew{t^{n+1}}
\def\tnhalf{t^{n+1/2}}

\def\vbar{\bar v}

\def\nuone{\nu+1}
\def\Inu{I_{\nu}}
\def\Inuhalf{I_{\nuhalf}}
\def\vnu{v_\nu}

\def\vnuhalfb{v}
\def\vnuhalf{\vnuhalfb_{\nuhalf}}
\def\vnup{v_\nu{}\!'}
\def\vnuonep{v{}'_{\nuone}}
\def\vnupp{v_\nu{}\!''}

\def\anub{a}
\def\anu{\anub_{\nu+1/2}}

\def\hatv{\widehat{v}}
\def\phik{\phi_{k}}

\def\eqspacea{\vspace{-0.1cm}}
\def\eqspaceb{\vspace{-0.08cm}}

\title{\bf  High Resolution Methods \vskip -2mm
for Time Dependent Problems \vskip -2mm with Piecewise Smooth Solutions\vskip 6mm}

\author{Eitan Tadmor\vspace*{-0.5cm}\thanks{University of California Los-Angeles, Department of Mathematics, Los Angeles, CA, USA 90095, and University of Maryland, Department of Mathematics, Center for Scientific Computation and Mathematical Modeling (CSCAMM) and Institute for Physical Science \& Technology (IPST), College Park, MD, USA 20742. E-mail: tadmor@math.umd.edu}}
\date{\vspace{-8mm}}

\begin{document}
\maketitle

\thispagestyle{first} \setcounter{page}{747}

\begin{abstract}\vskip 3mm
A trademark of nonlinear, time-dependent, convection-dominated problems is the spontaneous formation of non-smooth macro-scale features, like shock discontinuities and
non-differentiable kinks, which pose a challenge for high-resolution computations.
We overview recent developments of modern computational methods for
the approximate solution of such problems.
In these computations, one seeks piecewise smooth solutions which are realized by
finite dimensional projections. Computational methods in this context
can be classified into two main categories, of local and global methods.
Local methods are expressed in terms of point-values (--- Hamilton-Jacobi equations), cell averages (--- nonlinear conservation laws), or higher localized moments. Global methods are expressed in terms of global basis functions.\\

High resolution central schemes will be discussed as a prototype example for local methods.
The family of central schemes offers high-resolution ``black-box-solvers''
to an impressive range of such nonlinear problems. The main ingredients here are detection of spurious extreme values, non-oscillatory reconstruction in the directions of smoothness,
numerical dissipation and quadrature rules. Adaptive spectral viscosity will be discussed as an example for high-resolution global methods. The main ingredients here are detection of edges in spectral data, separation of scales, adaptive reconstruction, and spectral viscosity.

\vskip 4.5mm


\end{abstract}

\vskip 12mm

\section{Introduction --- piecewise smoothness} \label{section 1}\setzero
\vskip-5mm \hspace{5mm }

A trademark of nonlinear time-dependent convection-dominated problems is the spontaneous formation of non-smooth macro-scale features which challenge high-resolution computations.
A prototype example is the formation of shock discontinuities in nonlinear conservation laws,

\eqspacea
\be\label{eq:cons_law}
\frac{\part}{\part t} u(x,t)+\nabla_x\cdot f(u(x,t))=0, \quad u:=(u_1,\ldots,u_m)^\top.
\ee
\eqspaceb
\noindent
It is well known, e.g., \cite{Lax}, that solutions of
(\ref{eq:cons_law}) cease to be continuous, and (\ref{eq:cons_law})
should be interpreted in a weak sense with the derivatives on the left as Radon measures. This requires clarification. If $u(x,t)$ and $v(x,t)$ are two admissible solutions of (\ref{eq:cons_law}) then the following stability estimate is sought (here and below $\alpha, \beta,...$ stand for different generic constants),
\eqspacea
\be\label{eq:well_posed}
\|u(\cdot,t)-v(\cdot,t)\| \leq \alpha_t \|u(\cdot,0)-v(\cdot,0)\|.
\ee
\eqspaceb
Such estimates with different norms, $\|\cdot\|$, are playing a key role in
the linear setting --- both in theory and computations.
For linear hyperbolic systems, for example, (\ref{eq:well_posed})
is responsible for the usual $L^2$-stability theory, while
the stability of parabolic systems is often
measured by the $L^\infty$-norm, consult \cite{Richtmyer-Morton}.
But for nonlinear conservation laws, (\ref{eq:well_posed}) fails for any $L^p$-norm
with $p>1$. Indeed, comparing  $u(x,t)$ with any fixed translation
of it, $v(x,t):=u(x+h,t)$, the $L^p$ version of (\ref{eq:well_posed}) implies
\eqspacea
\[
\|\Delta_{+h}u(\cdot,t)\|_{L^p(\R^d)} \leq \alpha_t  \|\Delta_{+h}u(\cdot, 0)\|_{L^p(\R^d)},
\quad \Delta_{+h}u(\cdot,t):=u(\cdot+h,t)-u(\cdot,t).
\]
\eqspaceb
\noindent
For smooth initial data, however, the  bound on the right yields $\|\Delta_{+h}u(\cdot,t)\|_{L^p} \leq \alpha_t |h|$,
which in turn, for $p>1$, would lead to the contradiction that $u(\cdot,t)$ must remain continuous.
Therefore, conservation laws cannot satisfy the $L^p$-stability estimate (\ref{eq:well_posed})
after their finite breakdown time, except for the case $p=1$.
The latter leads to Bounded Variation (BV) solutions,
$\|u(\cdot,t)\|_{BV} :=
\sup \|\Delta_{+h}u(\cdot,t)\|_{L^1}/|h| \leq \alpha_t < \infty$,
whose derivatives are interpreted as the Radon measures mentioned above.
BV serves as the standard regularity space for admissible solutions of (\ref{eq:cons_law}).
A complete BV theory for scalar conservation laws, $m=1$, was developed Kru\v{z}kov.
Fundamental results on BV solutions of one-dimensional systems, $d=1$,
were obtained by P. Lax, J. Glimm, and others. Consult \cite{Bianchini-Bressan} for recent developments. Relatively little is known for general $(m-1)\times (d-1) >0$, but cf., \cite{Majda}.

We argue that the space of BV functions is still too large to describe the
approximate solutions of (\ref{eq:cons_law}) encountered in computations. Indeed, in such computations one does not 'faithfully' realize arbitrary BV functions but rather, piecewise smooth solutions. We demonstrate this point in the context of scalar approximate solutions,
$v^h(x,t)$, depending on a small computational scale $h\sim 1/N$. A typical error estimate for such approximations reads, \cite{CIME-lectures}
\eqspacea
\be\label{eq:BV_rate}
\|v^h(\cdot,t)-u(\cdot,t)\|_{L_{loc}^1(\R)} \leq \|v^h(\cdot,0)-u(\cdot,0)\|_{L_{loc}^1(\R)} + \alpha_t h^{1/2}.
\ee
\eqspaceb
The convergence rate of order $1/2$ is a well understood linear phenomena, which is observed
in computations\footnote{Bernstein polynomials, $B_N(u)$, provide a classical example of first-order monotone approximation with $L^1$-error of order $(\|u\|_{BV}/N)^{1/2}$. The general linear setting is similar, with improved rate $\sim h^{r/(r+1)}$ for $r$-order schemes.}.
The situation in the nonlinear case is different.
The optimal convergence rate for arbitrary BV initial data is still of order one-half, \cite{Sabac},
but actual computations exhibit higher-order convergence rate. The apparent difference between  theory and computations is resolved once we take into account piecewise smoothness. We can quantify piecewise smoothness in the simple scalar convex case, where the number of shock discontinuities of $u(\cdot,t)$ is bounded by the finitely many inflection points of the initial data, $u(x,0)$. In this case, the singular support of $u(\cdot,t)$ consists of
finitely many points where ${\cal S}(t)=\{x\ | \ \part_x u(x,t)\downarrow -\infty\}$.
Moreover, the solution in between those point discontinuities is as smooth as the initial data permit, \cite{Tadmor-Tassa}, namely
\eqspacea
\[
\sup_{x \in {\cal S}_L(t)} |\part_x^p u(x,t)| \leq e^{pLT}
\!\!\sup_{x \in {\cal S}_L(0)} |\part_x^pu(x,0)| + Const_L, \
{\cal S}_L(t):=\{x \ | \ \part_x u(x,t) \geq -L\}.
\]
\eqspaceb
If we let $d(x,t):=dist(x,{\cal S}(t))$ denote the distance to ${\cal S}(t)$,
then according to \cite{Tadmor-Tang}, the following  pointwise error estimate holds,
$|v^h(x,t)-u(x,t)| \leq \alpha_t  h/d(x,t)$, and integration yields the first-order convergence rate
\eqspacea
\be\label{eq:piecewise_rate}
\|v^h(\cdot,t)-u(\cdot,t)\|_{L_{loc}^1(\R)} \leq \alpha_t h|\log(h)|.
\ee
\eqspaceb
There is no contradiction between the optimality of (\ref{eq:BV_rate}) and (\ref{eq:piecewise_rate}). The former applies to arbitrary BV data, while the latter is restricted to piecewise smooth data and it is the one encountered in actual computations.
The general situation is of course, more complicated,  with a host of macro-scale features
which separate between regions of smoothness. Retaining the invariant properties of piecewise smoothness in general problems is a considerable challenge for high-resolution methods.

\section{A sense of direction} \label{section 2}
\setzero\vskip-5mm \hspace{5mm }

A computed approximation is a finite dimensional realization of
an underlying solution which, as we argue above, is viewed as
a piecewise smooth solution. To achieve higher accuracy, one should
extract more information from the smooth parts of the solution.  Macro-scale features of non-smoothness like shock discontinuities, are identified here as barriers for propagation of smoothness, and stencils which discretize (\ref{eq:cons_law}) while crossing discontinuities are excluded because of spurious Gibbs' oscillations. A high resolution scheme should sense
the direction of smoothness.

Another sense of directions is dictated by the propagation of information governed by  convective equations.
Discretizations of such equations fall into one of two, possibly overlapping categories.
One category of so-called upwind schemes consists of stencils which are fully aligned with the local direction of propagating waves. Another category of so-called central schemes consists
of two-sided stencils, tracing  both right-going and left-going waves.
A third possibility of stencils which discretize (\ref{eq:cons_law}) 'against the wind' is excluded because of their inherent instability, \cite{Richtmyer-Morton}. A stable scheme should sense the direction of propagation.

At this stage, high resolution stable schemes should compromise between two
different sets of directions, where propagation and smoothness might disagree.
This require essentially nonlinear schemes, with stencils which adapt their
sense of direction according to the computed data. We shall elaborate the details
in the context of high-resolution central scheme.

\section{Central schemes} \label{section 3} \setzero\vskip-5mm \hspace{5mm }

We start with a quotation from  \cite[\S12.15]{Richtmyer-Morton}, stating
``In 1959, Godunov described an ingenious method for one-dimensional problems with shocks''.
Godunov scheme is in the crossroads between the three major types of local discretizations, namely, finite-difference, finite-volume and finite-element methods.
The ingenuity of Godunov's approach, in our view, lies with the evolution of a
globally defined approximate solution, $v^h(x,\tn)$, replacing the prevailing approach at that time of an approximate solution which is realized by
its discrete gridvalues, $\vnu(\tn)$. This enables us to pre-process, to evolve and to post-process a globally defined approximation, $v^h(x,t^n)$. The main issue is how to 'manipulate' such piecewise smooth approximations while preserving the desired non-oscillatory invariants.

Godunov scheme was originally formulated in the context of nonlinear conservation laws, where an approximate solution is realized in terms of a
first-order accurate, piecewise-constant approximation
\eqspacea
\[
v^h(x,\tn):={\cal A}_{h} v(x,\tn):= \sum_\nu \vbar_\nu(\tn) 1_{\Inu}(x), \quad
\vbar_\nu(\tn):= \frac{1}{|\Inu|}\int_{\Inu} v(y,\tn)dy.
\]
\eqspaceb
The cell averages, $\vbar_{\nu}$, are evaluated over the equi-spaced cells,
$\Inu:=\{x \setminus |x-\xnu|\leq h/2 \}$ of uniform width $h\equiv \Delta x$.
More accurate Godunov-type schemes were devised using higher-order piecewise-polynomial projections. In the case of one-dimensional equi-spaced grid, such projections take the form
\eqspacea
\[
 {\cal P}_{h} v(x)= \sum_\nu p_\nu(x) 1_{\Inu}(x), \qquad
p_\nu(x)=\vnu + \vnup\Big(\frac{x-\xnu}{h}\Big) +
    \frac{1}{2}\vnupp\Big(\frac{x-\xnu}{h }\Big)^2+\ldots.
\]
\eqspaceb
Here, one pre-process the first-order cell averages in order to reconstruct accurate pointvalues, $\vnu$, and say, couple of numerical derivatives $\vnup/h, \vnupp/h^2$, while the original cell averages, $\{\vbar_\nu\}$,  should be preserved, ${\cal A}_h{\cal P}_hv^h={\cal A}_hv^h$.  The main issue is extracting information in the direction of smoothness. For a prototype example,
let $\Delta_+ \vbar_\nu$ and $\Delta_- \vbar_\nu$ denote the usual forward  and backward differences, $\Delta_\pm v_\nu:=\pm(\vbar_{\nu\pm 1}-\vbar_\nu)$. Starting with the given cell averages, $\{\vbar_\nu\}$, we set $\vnu=\vbar_\nu$, and compute
\eqspacea
\be\label{eq:minmod}
\vnup=mm(\Delta_+\vbar_\nu, \Delta_- \vbar_\nu), \quad
mm(z_1,z_2):=\frac{sgn(z_1)+sgn(z_2)}{2}\min\{|z_1|,|z_2|\}.
\ee
\eqspaceb
The resulting piecewise-linear approximation  is a second-order accurate, Total Variation Diminishing (TVD) projection, $\|{\cal P}_{h} v^h(x)\|_{BV} \leq  \|v^h(x)\|_{BV}$.
This recipe of so-called minmod numerical derivative, (\ref{eq:minmod}), is a representative for a large library of non-oscillatory, high-resolution limiters.
Such limiters dictate discrete stencils in the direction of smoothness and hence,
are inherently nonlinear. Similarly, nonlinear adaptive stencils are used in conjunction with higher-order methods. A description of the pioneering contributions in this direction by Boris \& Book, A. Harten, B. van-Leer and P. Roe  can be found in \cite{LeVeque}.
The advantage of dealing with  globally defined approximations is the ability to pre-process, to post-process and in particular, to evolve such approximations.  Let $u(x,t)=u^h(x,t)$ be the exact solution of (\ref{eq:cons_law}) subject to
$u^h(x,t^n)={\cal P}_hv^h(x,t^n)$. The exact solution lies of course outside the finite computational space, but it could be realized in terms of its exact cell averages, $v^h(x,\tnew)=\sum_\nu \vbar_\nu(\tnew)1_{\Inu}(x)$.
Averaging is viewed here a simple post-processing.
Two prototype examples are in order.\\
Integration of (\ref{eq:cons_law}) over control volume
$\Inu\times[\tn,\tnew]$, forms a local stencil which balances  between the new averages, $\{ \vbar_\nu(\tnew)\}$, the old ones, $\{\vbar_{\nu+k}(\tn)\}$,
and the fluxes across the interfaces along $x_{\nu\pm1/2}\times [\tn,\tnew]$.
In this case, the solution along these discontinuous interfaces is resolved in terms
of Riemann solvers. Since one employs here an exact evolution, the resulting Godunov-type schemes are upwind scheme. The original Godunov scheme based on piecewise constant projection
is the forerunner of all upwind schemes. As an alternative approach, one can realize the solution $u(x,\tnew)$, in terms
of its exact staggered averages, $\{\vbar_{\nuhalf}(\tnew)\}$.
Integration of (\ref{eq:cons_law}) over the control volume
$\Inuhalf\times[\tn,\tnew]$ subject to piecewise quadratic data given at $t=\tn$,
$u(x,\tn)={\cal P}_h v^h(x,\tn)$, yields
\eqspacea
\begin{eqnarray}
\vbar_{\nuhalf}(\tnew) & = & \frac{1}{2}\Big(\vbar_\nu(\tn)+\vbar_{\nu+1}(\tn)\Big) + \nonumber \\
    &  \ \ +  & \frac{1}{8}\Big(\vnup(\tn)-\vnuonep(\tn)\Big)
+\frac{\Delta t}{\Delta x}\Big(F^{n+1/2}_{\nuone}-F^{n+1/2}_\nu\Big), \label{eq:central_corrector}
\end{eqnarray}
\eqspaceb where $F^{n+1/2}_\nu$ stands for the averaged flux,
$F^{n+1/2}_\nu = \int_{\tn}^{\tnew} f(u(\xnu,\tau)d\tau/\Delta t$.
\linebreak
Thanks to the staggering of the grids, one encounters
smooth interfaces $\xnu\times [\tn,\tnew]$, and the intricate
(approximate) Riemann solvers are replaced by simpler quadrature
rules. For second-order accuracy, for example, we augment
(\ref{eq:central_corrector}) with the mid-point quadrature
\eqspacea \be\label{eq:2nd_central_predctor} F^{n+1/2}_\nu =
f(\vnu(\tnhalf)), \quad
    \vnu(\tnhalf)=\vnu(\tn)-\frac{\Delta t}{2\Delta x}f(\vnu(\tn))'.
\ee
\eqspaceb
Here, the prime on the right is understood in the usual sense of numerical differentiation of a gridfunction -- in this case the flux $\{f(\vnu(\tn))\}_\nu$. The resulting second-order central scheme (\ref{eq:2nd_central_predctor}),(\ref{eq:central_corrector}) was introduced in
\cite{Nessyahu-Tadmor}. It amounts to a simple predictor-corrector, non-oscillatory  high-resolution Godunov-type scheme. For systems, one implements numerical differentiation
for each component separately. Discontinuous edges are detected wherever cell-averages
form  new extreme values, so that $\vnup(\tn)$ and $f(\vnu(\tn))'$ vanish, and
(\ref{eq:2nd_central_predctor}),(\ref{eq:central_corrector}) is reduced to the
forerunner of all central schemes --- the celebrated first-order Lax-Friedrichs scheme,
\cite{LeVeque}. This first-order stencil is localized to the neighborhood of discontinuities,
and by assumption, there are finitely many them. In between those discontinuities, differentiation in the direction of smoothness
restores second-order accuracy. This retains the overall high-resolution of the scheme.
Consult Figure \ref{fig:electron_velocity} for example.

Similarly, higher-order quadrature rules can be used in connection with
higher-order projections, \cite{Bianco-Puppo-Russo}, \cite{Levy-Puppo-Russo}. A third-order simulation is presented in
Figure {\ref{fig:MHD}. Finite-volume and finite-element extensions
in  several space dimensions are realized over general, possibly unstructured control volumes, $\Omega_\nu\times [\tn,\tnew]$, which are adapted to handle general geometries.
Central schemes for 2D Cartesian grids were introduced in \cite{Jiang-Tadmor},
and extended to unstructured grids in \cite{Arminjon-Viallon}. A similar framework based on triangulated grids for high-resolution central approximations of Hamilton-Jacobi equations was described in \cite{CIME-lectures} and the references therein.

Central schemes enjoy the advantage of simplicity --  they are free of (approximate) Riemann solvers, they do not require dimensional splitting, and they apply to arbitrary flux
functions\footnote{An instructive example is provided by gasdynamics equations with tabulated equations of state.} without specific references to eigen-decompositions, Jacobians etc. In this context, central schemes offer a ``black-box solvers'' for an impressive variety of convection-dominated problems. At the same time, the central framework
maintain high-resolution by pre -and post-processing
in the direction of smoothness. References to diverse applications such as simulations of semi-conductors models, relaxation problems, geometrical optics and multiphase computations, incompressible flows, polydisperse suspensions, granular avalanches
 MHD equations and more  can be found at \cite{centralstation}.

\vspace{-0.4cm}
\renewcommand{\thefigure}{\arabic{figure}}
\begin{center}
\begin{minipage}[H]{5.8cm}
\begin{figure}[H]
\includegraphics[height=4.8cm,width=6.8cm]{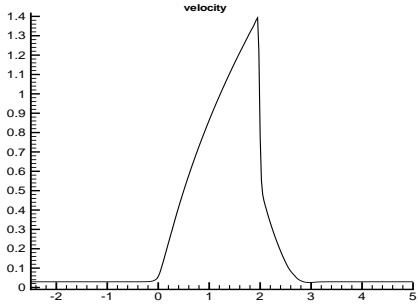}
\caption{Second-order central scheme simulation of semiconductor device governed by
1D Euler-Poisson equations. Electron velocity in $10^7$ cm/s with $N=400$ cells.}
\label{fig:electron_velocity}
\end{figure}
\end{minipage}\ \hspace{0.4cm}\
\renewcommand{\thefigure}{\arabic{figure}}
\begin{minipage}[H]{5.8cm}
\begin{figure}[H]
\includegraphics[height=4.7cm,width=5.8cm]{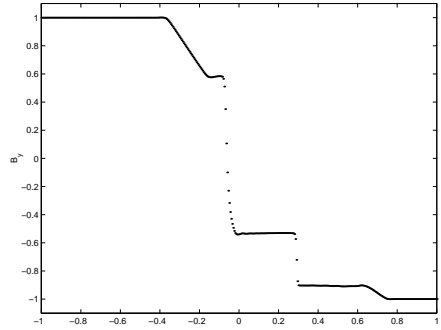}
\caption{Third-order central scheme simulation of 1D MHD Riemann problem.
with $N=400$ cells. The $y$-magnetic field at $t=0.2$.}
\label{fig:MHD}
\end{figure}
\end{minipage}
\end{center}

The numerical viscosity present in central schemes is of order ${\cal O}\frac{(\Delta x)^{2r}}{\Delta t}$. It is suitable for the convective regime where $\Delta t \sim \Delta x$,
but it is  excessive when a small time step is enforced, e.g., due to the presence of diffusion terms. To overcome this difficulty, a new family
of central schemes was introduced in \cite{Kurganov-Tadmor} and was further refined in \cite{Kurganov-Noelle-Petrova}. Here,  the previous staggered control volumes,
$I_{\nuhalf}\times [\tn,\tnew]$,  is replaced by the smaller --- and hence less dissipative,
$J_{\nuhalf} \times [\tn,\tnew]$, where
$J_{\nuhalf}:=\xnuhalf+\Delta t\times[\anub^-, \anub^+]$ encloses the maximal cone
of propagation, $\anub^\pm\equiv\anu^\pm= {{\rm max}\choose {\rm min}}_k \lambda_k^\pm(f_u)$.
The fact that the staggered grids are ${\cal O}(\Delta t)$ away from each other,
yields central stencils with numerical viscosity of order ${\cal O}(\Delta x^{2r-1})$.
Being independent of $\Delta t$ enables us to pass to the limit
$\Delta t \downarrow 0$. The resulting semi-discrete high-resolution central scheme reads
$\dot{\vbar}_\nu(t) =  -(f_{\nuhalf}(t)-f_{\numhalf}(t))/\Delta x $,
with a numerical flux, $f_{\nuhalf}$, expressed in terms of the  reconstructed pointvalues,
$\vnuhalfb^\pm\equiv \vnuhalf^\pm= {\cal P}_h v^h(\xnuhalf\pm, t)$,
\eqspacea
\be\label{eq:semi_discrete_flux}
 f_{\nuhalf}(t):= \frac{\anub^+ f(\vnuhalfb^-)-\anub^- f(\vnuhalfb^+)}{\anub^+ - \anub^-}
 + \anub^+\anub^-\frac{\vnuhalfb^+ - \vnuhalfb^-}{\anub^+-\anub^-}.
\ee
\eqspaceb
Instructive examples are provided in Figures \ref{fig:convection-diffusion},
\ref{fig:2D_Riemann}.

\vspace{-0.4cm}
\renewcommand{\thefigure}{\arabic{figure}}
\begin{center}
\begin{minipage}[H]{5.8cm}
\begin{figure}[H]
\includegraphics[height=4.0cm,width=5.8cm]{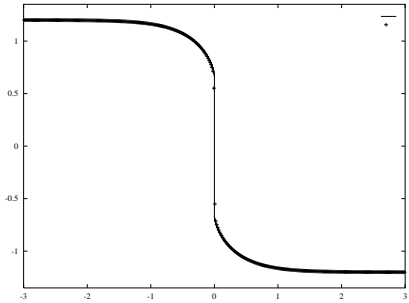}
\caption{Convection-diffusion eq.
$u_t + (u^2)_x/2= (u_x/\sqrt{1+u_x^2})_x$ simulated by
(\ref{eq:semi_discrete_flux}),(\ref{eq:minmod}), with $400$ cells.}
\label{fig:convection-diffusion}
\end{figure}
\end{minipage}\ \hspace{0.4cm}\
\renewcommand{\thefigure}{\arabic{figure}}
\begin{minipage}[H]{5.9cm}
\begin{figure}[H]
\includegraphics[height=4.0cm,width=5.8cm]{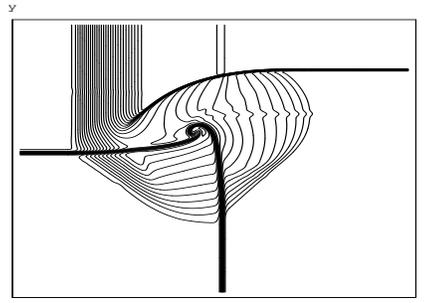}
\caption{Third-order  central scheme for 2D Riemann problem.
Density contour lines with $400\times 400$ cells.}
\label{fig:2D_Riemann}
\end{figure}
\end{minipage}
\end{center}

\section{Adaptive spectral viscosity methods} \label{section 4} \setzero\vskip-5mm \hspace{5mm }

Godunov-type methods are based on zeroth-order moments of (\ref{eq:cons_law}). In each time step, one evolves one piece of information per spatial cell --- the cell average. Higher accuracy is restored by numerical differentiation in the direction of smoothness. An alternative approach is to compute higher-order moments, where the cell averages, $\vbar_\nu$, and say, couple of numerical derivatives, $\vnup, \vnupp$ are evolved in time. Prototype examples include  discontinuous Galerkin and streamline-diffusion methods, where several local moments
per computational cell are evolved in time, consult \cite{CIME-lectures} and the references therein. As the number of projected moments increase, so does the size of the computational stencil. At the limit, one arrives at spectral methods based on global projections,
\eqspacea
\[
v_N(x,t)= {\cal P}_N v(x,t):= \sum_{|k|\leq N} \hatv_k(t)\phik(x), \quad \hatv_k:=\langle v(\cdot,t),\phik(\cdot)\rangle.
\]
\eqspaceb
Here, $\phi_k(x)$ are global basis functions, $\phik=e^{ikx}, T_k(x), ...$
and $\hatv_k$ are the moments induced by the appropriately weighted, possibly discrete  inner-product $\langle\cdot,\cdot\rangle$. Such global projections enjoy superior resolution --- the error $\|v(\cdot)-{\cal P}_Nv(\cdot)\|$ decays as fast as the global smoothness of $v(\cdot)$ permits. With piecewise smooth solutions, however, we encounter first-order spurious Gibbs' oscillations throughout the computational domain. As before, we should be able to pre- and post-process piecewise smooth projections, recovering their high accuracy in the direction of smoothness. To this end, local limiters are replaced by global concentration kernels of the form, \cite{Gelb-Tadmor},
$K^\eta_Nv_N(x) = \frac{\pi \omega(x)}{N}\sum_{|k|\leq N} \eta\Big(\frac{|k|}{N}\Big)\hatv_k (\phik(x))_x$, where $\eta(\cdot)$ is an arbitrary unit mass $C^\infty[0,1]$ function at our disposal.
Detection of edges is facilitated by separation between the ${\cal O}(1)$ scale in the neighborhoods of edges and ${\cal O}(h^r)$ scales in regions of smoothness, $K_N v_N(x) \sim [v(x)] + {\cal O}(Nd(x))^{-r}$. Here, $[v(x)]$ denotes the amplitude of the jump discontinuity at $x$ (-- with vanishing amplitudes signaling smoothness), and $d(x)$ is the distance to the singular support of $v(\cdot)$. Two prototype examples in the $2\pi$-periodic setup are in order.
With $\eta(\xi)\equiv 1$ one recovers a first-order concentration kernel due to Fejer. In \cite{Gelb-Tadmor} we introduced  the concentration kernel, $\eta^{exp}(\xi) :\sim exp\{\beta/\xi(1-\xi)\}$,  with exponentially fast decay into regions of smoothness.
Performing the minmod limitation,
(\ref{eq:minmod}), $mm\{K_N^1 v_N(x), K_N^{exp}v_N(x)\}$, yields an adaptive, essentially non-oscillatory  edge detector with enhanced separation of scales  near jump discontinuities.
Once macro-scale features of non-smoothness were located, we turn to reconstruct the
information in the direction of smoothness. This could be carried out either in the physical space using adaptive mollifiers, $\Psi_{\theta,p}$, or by nonlinear adaptive filters, $\sigma_{p}$. In the $2\pi$ periodic case, for example, we set
$ \Psi*{\cal P}_Nv(x):= \langle\Psi_{\theta,p}(x-\cdot), {\cal P}_Nv(\cdot)\rangle$
where $\Psi$ is expressed in terms of the Dirichlet kernel,
$D_p(y):=\frac{\sin(p+1/2)\pi y}{2\sin(\pi y/2)}$,
\eqspacea
\be\label{eq:adaptive_mollifier}
\Psi_{\theta,p}(y)= \frac{1}{\theta}\rho\Big(\frac{y}{\theta}\Big)D_p(\Big(\frac{y}{\theta}\Big), \quad
\rho\equiv \rho_\beta(y):=e^{\beta y^2/(y^2-1)}1_{[-1,1]}(y).
\ee
\eqspaceb
\noindent
Mollifiers encountered in applications maintain their finite accuracy by localization,
letting $\theta\downarrow 0$. Here, however, we seek superior accuracy by the process of cancellation with increasing $p\uparrow$. To guarantee that the reconstruction is supported in the direction of smoothness, we maximize $\theta$ in terms of the distance function, $\theta(x)=d(x)$,  so that we avoid crossing
discontinuous barriers. Superior accuracy is achieved by the adaptive choice $p :\sim d(x) N$, which yields, \cite{Tadmor-Tanner},
\eqspacea
\[
|\Psi_{\{d(x), d(x)N\}}*{\cal P}_Nv(x) - v(x)| \leq
Const. (d(x)N)^r e^{-\gamma \sqrt{d(x)N}}.
\]
\eqspaceb
The remarkable exponential recovery is due to the Gevrey regularity of
$\rho_\beta \in {\cal G}_2$ and Figure \ref{fig:adaptive_mollifier} demonstrates such
adaptive recovery with exponential convergence rate recorded in Figure \ref{fig:expo_error}. An analogous filtering procedure can be carried out in the dual Fourier space, and as before, it hinges on a filter with an adaptive degree $p$
\eqspacea
\[
\Psi*{\cal P}_Nv(x):=\sum_{|k|\leq N} \sigma_{p}\Big(\frac{|k|}{N}\Big)\hatv_k e^{ikx},
\quad \sigma_{p}(\xi)=e^{\beta\xi^p/(\xi^2-1)}, \ \ p\sim (d(x)N)^{r/r+1}.
\]
\eqspaceb

\vspace{-0.4cm}
\renewcommand{\thefigure}{\arabic{figure}}
\begin{center}
\begin{minipage}[H]{5.8cm}
\begin{figure}[H]
\includegraphics[height=4.0cm,width=5.8cm]{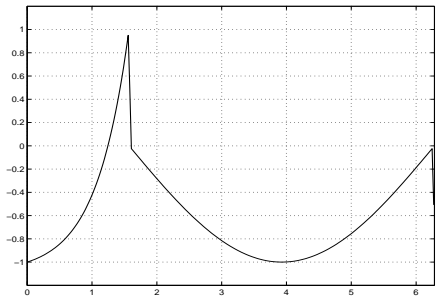}
\caption{Adaptive reconstruction of the piecewise smooth data from its using $N=128$-modes
using (\ref{eq:adaptive_mollifier}).}
\label{fig:adaptive_mollifier}
\end{figure}
\end{minipage}\ \hspace{0.4cm}\
\renewcommand{\thefigure}{\arabic{figure}}
\begin{minipage}[H]{5.8cm}
\begin{figure}[H]
\includegraphics[height=4.0cm,width=5.8cm]{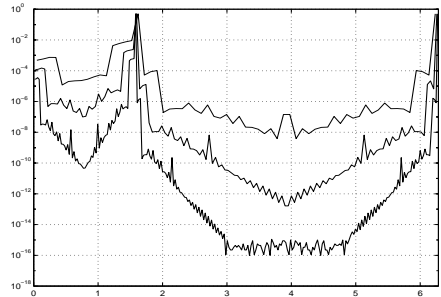}
\caption{Log error for an adaptive spectral mollifier
based on  $N=32,\ 64,\ \mbox{and} \ 128$ modes.}
\label{fig:expo_error}
\end{figure}
\end{minipage}
\end{center}

Equipped with this toolkit to process spectral projections of piecewise smooth data, we turn to consider their time evolution. Godunov-type methods are based evolution of cell averages.
Cell averaging is dissipative, but projections onto higher-order moments are not. The following example, taken from \cite{Tadmor-SV} shows what go wrong with global projections.
We consider the Fourier method for $2\pi$-periodic inviscid Burgers' equation,
$\part_t v_N(x,t)+ \part_x S_N(v_N^2(\cdot,t))/2=0$. Orthogonality implies that $\|v_N(\cdot,t)\|^2_{L^2}$ is conserved in time, and in particular, that the Fourier method  admits a weak limit, $v_N(\cdot,t) \rightharpoonup \vbar(t)$. At the same time, $\vbar(t)$ is not a strong limit, for otherwise it will contradict the strict  entropy dissipation associated with shock discontinuities. Lack of strong convergence indicates the persistence of
spurious dispersive oscillations, which is due to lack of entropy dissipation.
With this in mind, we turn to discuss the Spectral Viscosity (SV) method, as a
framework to stabilize the evolution of global projections without sacrificing their
spectral accuracy. To this end one augments the usual Galerkin procedure with high frequency
viscosity regularization
\eqspacea
\be\label{eq:SV}
\frac{\part}{\part t}v_N(x,t) +  {\cal P}_N\nabla_x\cdot f(v_N(\cdot,t))=
-N \sum_{|k|\leq N} \sigma\Big(\frac{|k|}{N}\Big)\hatv_k(t)\phi_k(x),
\ee
\eqspaceb
where $\sigma(\cdot)$ is a low-pass filter satisfying
$\sigma(\xi) \geq \big(|\xi|^{2s}-\frac{\beta}{N}\Big)^+$.
Observe that the SV is only activated on the highest portion of the spectrum,
with wavenumbers $|k|> \gamma N^{(2s-1)/2s}$. Thus, the SV method can be viewed as
a compromise between the stable
viscosity approximations corresponding to $s=0$ but restricted to first order,
and the spectrally accurate yet unstable spectral method corresponding to $s=\infty$.

\vspace{-0.6cm}
\renewcommand{\thefigure}{\arabic{figure}}
\begin{center}
\begin{minipage}[H]{5.8cm}
\begin{figure}[H]
\includegraphics[height=4.0cm,width=5.8cm]{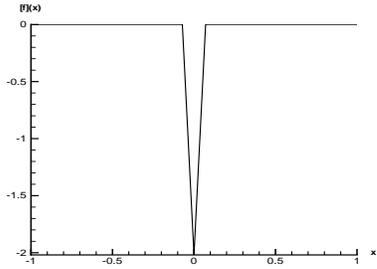}
\caption{Enhanced detection of edges with ~$v(x)$ given by $v(x)=-sgn x \cdot\cos(x+x\cdot sgnx/2)1_{[-\pi,\pi]}(x)$.}
\label{fig:detection}
\end{figure}
\end{minipage}\ \hspace{0.4cm}\
\renewcommand{\thefigure}{\arabic{figure}}
\begin{minipage}[H]{5.8cm}
\begin{figure}[H]
\includegraphics[height=4.0cm,width=5.8cm]{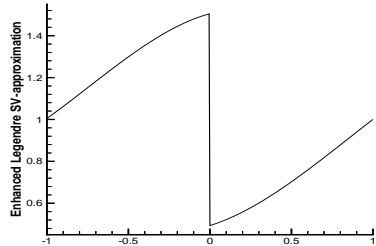}
\caption{Legendre SV solution of inviscid Burgers' equation.
Reconstruction in the direction smoothness.}
\label{fig:SV-method}
\end{figure}
\end{minipage}
\end{center}

\noindent
The additional SV on the right of (\ref{eq:SV}) is small enough to
retain the formal spectral accuracy of the underlying spectral approximation.
At the same time, SV is large  enough to enforce a sufficient
amount of entropy dissipation and hence, to prevent the unstable spurious Gibbs' oscillations. The original approach was introduced in \cite{Tadmor-SV} in conjunction with second-degree dissipation, $s=1$. Extensions to
several space variables, non-periodic expansions, further developments of hyper SV methods with $0 < s < \infty$, and various applications can be found in \cite{SVstation}. We conclude with an implementation of adaptive SV method for simple inviscid Burgers' in Figure \ref{fig:SV-method}.

\label{lastpage}

\end{document}